\documentclass[12pt,twoside, epsf]{article}

\usepackage{color,graphicx,times}

\makeatother

\let \ttorg \tt \def \tt{\ttorg \obeyspaces}

\begin{document}

\newcommand{\Across}{\raisebox{-0.25\height}{\includegraphics[width=0.5cm]{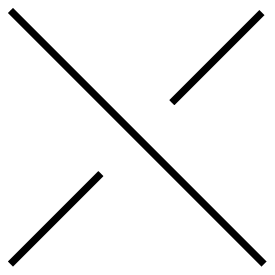}}}
\newcommand{\Asmooth}{\raisebox{-0.25\height}{\includegraphics[width=0.5cm]{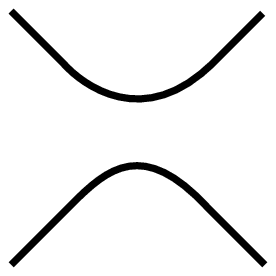}}}
\newcommand{\Bsmooth}{\raisebox{-0.25\height}{\includegraphics[width=0.5cm]{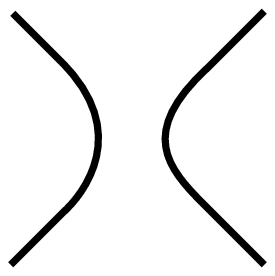}}}
\newcommand{\Rcurl}{\raisebox{-0.25\height}{\includegraphics[width=0.5cm]{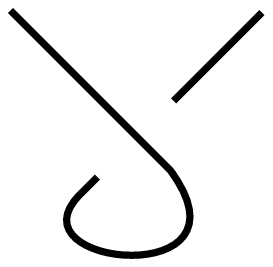}}}
\newcommand{\Lcurl}{\raisebox{-0.25\height}{\includegraphics[width=0.5cm]{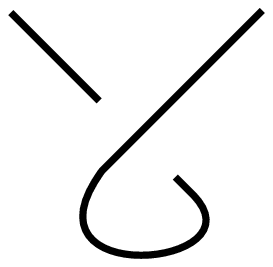}}}
\newcommand{\Arc}{\raisebox{-0.25\height}{\includegraphics[width=0.5cm]{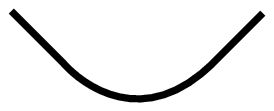}}}

\date{}

\title{\Large\bf Topological Quantum Information, Khovanov Homology and the Jones Polynomial}
\author{Louis
H. Kauffman\\ Department of Mathematics, Statistics \\ and Computer Science (m/c
249)    \\ 851 South Morgan Street   \\ University of Illinois at Chicago\\
Chicago, Illinois 60607-7045\\ $<$kauffman@uic.edu$>$\\}

\maketitle

\thispagestyle{empty}

\subsection*{\centering Abstract}

{\em This paper is dedicated to Anatoly Libgober on his 60-th birthday. }

\section{Introduction}
In this paper we give a quantum statistical interpretation for the bracket polynomial state sum $\langle K \rangle $ and for the
Jones polynomial $V_{K}(t).$ We use this quantum mechanical interpretation to give a new quantum algorithm for computing the Jones polynomial.
This algorithm is useful for its conceptual simplicity  and it applies to all values of the polynomial variable that lie on the unit circle in the 
complex plane. Letting
${\cal C}(K)$ denote the Hilbert space for this model, there is a natural unitary transformation
$$U: {\cal C}(K) \longrightarrow {\cal C}(K)$$ such that $$\langle K \rangle = <\psi|U|\psi>$$ where $| \psi \rangle$ is a sum over basis states
for ${\cal C}(K).$ The quantum algorithm comes directly from this formula via the Hadamard Test.
We then show that the framework for our quantum model for the bracket polynomial
is a natural setting for Khovanov homology. The Hilbert space ${\cal C}(K)$ of our model has basis in one-to-one correspondence with 
the enhanced states of the bracket state summmation and is isomorphic with the chain complex for Khovanov homology with coefficients in the 
complex numbers. We show that for the Khovanov boundary operator $\partial: {\cal C}(K) \longrightarrow {\cal C}(K),$ we have the relationship
$\partial U + U \partial = 0.$ Consequently, the operator $U$ acts on the Khovanov homology, and we therefore obtain a direct relationship between
Khovanov homology and this quantum algorithm for the Jones polynomial. The quantum algorithm given here is inefficient, and so it remains an open problem
to determine better quantum algorithms that involve both the Jones polynomial and the Khovanov homology.
\bigbreak

The paper is organized as follows. Section 2 reviews the structure of the bracket polynomial, its state summation, the use of enhanced states and the
relationship with the Jones polynomial. Section 3 describes the quantum statistical model for the bracket polynomial and refers to Section 7 (the Appendix)
for the details of the Hadamard Test and the structure of the quantum algorithm. Section 4 describes the relationship of the quantum model with Khovanov
homology axiomatically, without  using the specific details of the Khovanov chain complex that give these axioms. In Section 5 we construct the Khovanov
chain complex in detail and show  how some of its properties follow uniquely from the axioms of the previous section. In Section 6 we discuss how the
framework of this paper can be  generalized to other situations where the Hilbert space of a quantum information system is also a chain complex. We give an
example using DeRahm cohomology. In Section 7 (the Appendix) we detail the Hadamard test.
\bigbreak

\section{Bracket Polynomial and Jones Polynomial}
The bracket polynomial \cite{KaB} model for the Jones polynomial \cite{JO,JO1,JO2,Witten} is usually described by the expansion
$$\langle \Across \rangle=A \langle \Asmooth \rangle + A^{-1}\langle \Bsmooth \rangle$$
Here the small diagrams indicate parts of otherwise identical larger knot or link diagrams. The two types of smoothing (local diagram with no crossing) in
this formula are said to be of type $A$ ($A$ above) and type $B$ ($A^{-1}$ above). 

$$\langle \bigcirc \rangle = -A^{2} -A^{-2}$$
$$\langle K \, \bigcirc \rangle=(-A^{2} -A^{-2}) \langle K \rangle $$
$$\langle \Rcurl \rangle=(-A^{3}) \langle \Arc \rangle $$
$$\langle \Lcurl \rangle=(-A^{-3}) \langle \Arc \rangle $$
One uses these equations to normalize the invariant and make a model of the Jones polynomial.
In the normalized version we define $$f_{K}(A) = (-A^{3})^{-wr(K)} \langle K \rangle / \langle \bigcirc \rangle $$ 
where the writhe $wr(K)$ is the sum of the oriented crossing signs for a choice of orientation of the link $K.$ Since we shall not use oriented links
in this paper, we refer the reader to \cite{KaB} for the details about the writhe. One then has that $f_{K}(A)$ is invariant under the Reidemeister moves
(again see \cite{KaB}) and the original Jones polynonmial $V_{K}(t)$ is given by the formula $$V_{K}(t) = f_{K}(t^{-1/4}).$$ The Jones polynomial
has been of great interest since its discovery in 1983 due to its relationships with statistical mechanics, due to its ability to often detect the
difference between a knot and its mirror image and due to the many open problems and relationships of this invariant with other aspects of low 
dimensional topology.  
\bigbreak

\noindent {\bf The State Summation.} In order to obtain a closed formula for the bracket, we now describe it as a state summation.
Let $K$ be any unoriented link diagram. Define a {\em state}, $S$, of $K$  to be the collection of planar loops resulting from  a choice of
smoothing for each  crossing of $K.$ There are two choices ($A$ and $B$) for smoothing a given  crossing, and
thus there are $2^{c(K)}$ states of a diagram with $c(K)$ crossings.
In a state we label each smoothing with $A$ or $A^{-1}$ according to the convention 
indicated by the expansion formula for the bracket. These labels are the  {\em vertex weights} of the state.
There are two evaluations related to a state. The first is the product of the vertex weights,
denoted $\langle K|S \rangle .$
The second is the number of loops in the state $S$, denoted  $||S||.$
  
\noindent Define the {\em state summation}, $\langle K \rangle $, by the formula 

$$\langle K \rangle  \, = \sum_{S} <K|S> \delta^{||S||}$$
where $\delta = -A^{2} - A^{-2}.$
This is the state expansion of the bracket. It is possible to rewrite this expansion in other ways. For our purposes in 
this paper it is more convenient to think of the loop evaluation as a sum of {\it two} loop evaluations, one giving $-A^{2}$ and one giving 
$-A^{-2}.$ This can be accomplished by letting each state curve carry an extra label of $+1$ or $-1.$ We describe these {\it enhanced states}
below.
\bigbreak

\noindent {\bf Changing Variables.} Letting $c(K)$ denote the number of crossings in the diagram $K,$ if we replace $\langle K
\rangle$ by 
$A^{-c(K)} \langle K \rangle,$ and then replace $A^2$ by $-q^{-1},$ the bracket is then rewritten in the
following form:
$$\langle \Across \rangle=\langle \Asmooth \rangle-q\langle \Bsmooth \rangle $$
with $\langle \bigcirc\rangle=(q+q^{-1})$.
It is useful to use this form of the bracket state sum
for the sake of the grading in the Khovanov homology (to be described below). We shall
continue to refer to the smoothings labeled $q$ (or $A^{-1}$ in the
original bracket formulation) as {\it $B$-smoothings}. 
\bigbreak

\noindent {\bf Using Enhanced States.}
We now use the convention of {\it enhanced
states} where an enhanced state has a label of $1$ or $-1$ on each of
its component loops. We then regard the value of the loop $q + q^{-1}$ as
the sum of the value of a circle labeled with a $1$ (the value is
$q$) added to the value of a circle labeled with an $-1$ (the value
is $q^{-1}).$ We could have chosen the less neutral labels of $+1$ and $X$ so that
$$q^{+1} \Longleftrightarrow +1 \Longleftrightarrow 1$$
and
$$q^{-1} \Longleftrightarrow -1 \Longleftrightarrow X,$$
since an algebra involving $1$ and $X$ naturally appears later in relation to Khovanov homology. It does no harm to take this form of labeling from the
beginning.
\bigbreak

Consider the form of the expansion of this version of the 
bracket polynonmial in enhanced states. We have the formula as a sum over enhanced states $s:$
$$\langle K \rangle = \sum_{s} (-1)^{i(s)} q^{j(s)} $$
where $i(s)$ is the number of $B$-type smoothings in $s$ and $j(s) = i(s) + \lambda(s)$, with $\lambda(s)$ equal to the number of loops labeled $1$ minus
the number of loops labeled $-1$ in the enhanced state $s.$
\bigbreak

One advantage of the expression of the bracket polynomial via enhanced states is that it is now a sum of monomials. We shall make use of this property
throughout the rest of the paper.
\bigbreak

\section{Quantum Statistics and the Jones Polynomial}
We now use the enhanced state summation for the bracket polynomial with variable $q$ to give a quantum formulation of the state sum.
{\it Let $q$ be on the unit circle in the complex plane.} (This is equivalent to letting the original bracket variable $A$ be on the unit
circle and equivalent to letting the Jones polynmial variable $t$ be on the unit circle.) Let ${\cal C}(K)$ denote the complex vector space
with orthonormal basis $ \{ |s\rangle$ \} where $s$ runs over the enhanced states of the diagram $K.$ The vector space ${\cal C}(K)$
is the (finite dimensional) Hilbert space for our quantum formulation of the Jones polynomial. 
While it is customary for a Hilbert space to be written with the letter
$H,$ we do not follow that convention here, due to the fact that we shall soon regard ${\cal C}(K)$ as a chain complex and take its
homology. One can hardly avoid using  $\cal{H}$ for homology.
\bigbreak

\noindent With $q$ on the unit circle, we define a unitary transformation $$U: {\cal C}(K) \longrightarrow {\cal C}(K)$$ by
the formula
$$U |s \rangle = (-1)^{i(s)} q^{j(s)} |s \rangle$$ for each enhanced state $s.$ Here $i(s)$ and $j(s)$ are as defined in the previous 
section of this paper.
\bigbreak

\noindent Let $$| \psi \rangle = \sum_{s} |s \rangle.$$ The state vector  $ | \psi \rangle$ is the sum over the basis states of our Hilbert space ${\cal
C}(K).$ For convenience, we do not normalize $|\psi \rangle$ to length one in the Hilbert space ${\cal C}(K).$ 
We then have the 
\bigbreak

\noindent {\bf Lemma.} The evaluation of the bracket polynomial is given by the following formula
$$\langle K \rangle = \langle \psi | U | \psi \rangle .$$ 
\bigbreak

\noindent {\bf Proof.} $$\langle \psi | U | \psi \rangle = \sum_{s'} \sum _{s} \langle s'| (-1)^{i(s)}q^{j(s)} |s \rangle = \sum_{s'} \sum
_{s}(-1)^{i(s)}q^{j(s)} \langle s'| s \rangle $$
$$ =  \sum_{s}(-1)^{i(s)}q^{j(s)} = \langle K \rangle,$$ since $$\langle s'| s \rangle = \delta(s,s')$$
where $\delta(s,s')$ is the Kronecker delta, equal to $1$ when $s = s'$ and equal to $0$ otherwise. //
\bigbreak

Thus the bracket polyomial evaluation is a quantum amplitude for the measurement of the state
$U | \psi \rangle$ in the $\langle \psi |$ direction. Since $\langle \psi | U | \psi \rangle$ can be regarded as a diagonal element of the 
transformation $U$ with respect to a basis containing $|\psi \rangle,$ this formula can be taken as the foundation for a quantum algorithm that computes
the bracket of $K$ via the Hadamard test. See the Appendix to this paper, for a discussion of the Hadamard test and the corresponding quantum algorithm.
\bigbreak

A few words about quantum algorithms are appropropiate at this point. In general one begins with an intial state $| \psi \rangle$ and a unitary
transformation $U.$ Quantum processes are modeled by unitary transformations, and it is in principle possible to create a physical process
corresponding to any given unitary transformation. In practice, one is limited by the dimensions of the spaces involved and by the fact that 
physical quantum states are very delicate and subject to decoherence. Nevertheless, one designs
quantum algorithms at first by finding unitary operators that represent the informaton in the problem one wishes to calculate. Then the quantum part
of the quantum computation is the physical process that produces the state  $U | \psi \rangle$ from the initial state  $| \psi \rangle.$ At this point
$U | \psi \rangle$ is available for measurement. Measurement means interaction with the enviroment, and will happen in any case, but one intends a 
controlled circumstance in which the measurement can occur. If $\{ |e_{1} \rangle,\cdots |e_{n} \rangle \}$ is a basis for the Hilbert space
(here denoted by ${\cal C}(K)$) then we would have $$U| \psi \rangle = \sum_{k} z_{k} |e_{k} \rangle,$$ a linear combination of the basis elements.
On measuring with respect to this basis, each $|e_{k} \rangle$ corresponds to an observable outcome and {\it this outcome will occur with 
frequency $|z_{k}|^{2}/(\sum_{i}|z_{i}|^{2}).$} The form of quantum computation consists in repeatedly  preparing and running the process to find the
frequency with  which certain key observations occur.  In the case of our algorithm we are interested in the frequency of
measuring 
$| \psi \rangle$ itself and this corresponds to the inner product $\langle \psi |U | \psi \rangle.$ However, a direct measurement from $U | \psi \rangle$
will yield only approximations to the absolute square of $\langle \psi|U |\psi \rangle.$ For this reason we use the more refined Hadamard test as 
described in the Appendix. There is much more to be said about quantum algorithms in general and about quantum algorithms for the Jones polynomial in
particular. The reader can examine \cite{AJL,QCJP,3Strand,QKnots,NMR} for more information.
\bigbreak

It is useful to formalize the bracket evaluation as a quantum amplitude. This is a direct way to give a physical interpretation of the
bracket state sum and the Jones polynomial. Just how this process can be implemented physically depends upon the interpretation of the Hilbert
space ${\cal C}(K).$ It is common practice in theorizing about quantum computing and quantum information to define a Hilbert space in terms
of some mathematically convenient basis (such as the enhanced states of the knot or link diagram $K$) and leave open the possibility of a realization 
of the space and the quantum evolution operators that have been defined upon it. In 
principle any finite dimensional unitary operator can be realized by some physical system. In practice, this is the problem of constructing quantum 
computers and quantum information devices. It is not so easy to construct what can be done in principle, and the quantum states that are produced
may be all too short-lived to produce reliable computation. Nevertheless, one has the freedom to create spaces and operators on the mathematical level
and to conceptualize these in a quantum mechanical framework. The resulting structures may be realized in nature and in present or future
technology. In the case of our Hilbert space associated with the bracket state sum and its corresponding unitary transformation $U,$ there is rich extra 
structure related to Khovanov homology that we discuss in the next section. One hopes that in a (future) realization of these spaces and operators,
the Khovanov homology will play a key role in quantum information related to the knot or link $K.$
\bigbreak

There are a number of conclusions that we can draw from the formula $\langle K \rangle = \langle \psi | U | \psi \rangle .$ First of all, this formulation
constitutes a quantum algorithm for the computation of the bracket polynomial (and hence the Jones polynomial) at any specialization where the variable is
on the unit circle. We have defined a unitary transformation
$U$ and then shown that the bracket is an evaluation in the form $ \langle \psi| U | \psi \rangle.$ This evaluation can be computed via the Hadamard test
\cite{NC} and this gives the desired quantum algorithm. Once the unitary transformation is given as a physical construction, the algorithm will be as
efficient as any application of the Hadamard test. The present algorithm requires an exponentially increasing complexity of construction 
for the associated unitary transformation, since the dimension of the Hilbert space is equal to the $2^{e(K)}$ where $e(K)$ is the number of enhanced
states of the diagram $K$. (Note that $e(K) = \sum_{S} 2^{||S||}$ where $S$ runs over the $2^{c(K)}$ standard bracket states, $c(K)$ is the
number of crossings in the diagram and $||S||$ is the number of loops in the state $S.$) Nevertheless, it is significant that the Jones polynomial can be
formulated in such a direct way in terms of a quantum algorithm. By the same token, we can take the basic result of Khovanov homology that says that the
bracket is a graded Euler characteristic of the Khovanov homology as telling us that we are taking a step in the direction of a quantum algorithm for the
Khovanov homology itself. This will be discussed below.
\bigbreak

\section{Khovanov Homology and a Quantum Model for the Jones Polynomial}
In this section we outline how the Khovanov homology is related with our quantum model. This can be done essentially axiomatically, without
giving the details of the Khovanov construction. We give these details in the next section. The outline is as follows:
\begin{enumerate}
\item There is a boundary operator $\partial$ defined on the Hilbert space of enhanced states of a link diagram $K$
$$\partial:{\cal C}(K) \longrightarrow {\cal C}(K)$$ such that $\partial \partial = 0$ and so that if
${\cal C}^{i,j} = {\cal C}^{i,j}(K)$ denotes the subspace of ${\cal C}(K)$ spanned by enhanced states $| s \rangle$ with $i = i(s)$ and $j = j(s),$ then
$$\partial : {\cal C}^{ij} \longrightarrow {\cal C}^{i+1,j}.$$ That is, we have the formulas $$i(\partial |s \rangle) = i(| s \rangle) + 1$$ and
$$j(\partial | s \rangle ) = j(| s \rangle )$$ for each enhanced state $s.$ In the next section, we shall explain how the boundary operator is constructed.

\item {\bf Lemma. } By defining $U:{\cal C}(K) \longrightarrow {\cal C}(K)$ as in the previous section, via $$U|s \rangle = (-1)^{i(s)} q^{j(s)} |s
\rangle,$$ we have the following basic relationship between $U$ and the boundary operator $\partial:$ $$ U \partial + \partial  U = 0.$$
\smallbreak
 
\noindent {\bf Proof. } This follows at once from the definition of $U$ and the fact that $\partial$ preserves $j$ and increases $i$ to $i+1.$ //

\item  From this Lemma
we conclude that the operator $U$ acts on the homology of ${\cal C}(K).$ We can regard $H({\cal C}(K)) = Ker(\partial)/Image(\partial)$ as a new Hilbert
space on which the unitary  operator $U$ acts. In this way, the Khovanov homology and its relationship with the Jones polynomial has a natural quantum
context.

\item For a fixed value of $j$, $${\cal C}^{\bullet,j} = \oplus_{i} {\cal C}^{i,j}$$ is a subcomplex of ${\cal C}(K)$
with the boundary operator $\partial.$ Consequently, we can speak of the homology $H({\cal C}^{\bullet,j}).$ Note that the dimension of ${\cal C}^{ij}$ is
equal to the number of enhanced states $|s \rangle$ with $i = i(s)$ and $j = j(s).$ Consequently, we have
$$\langle K \rangle = \sum_{s} q^{j(s)}(-1)^{i(s)} = \sum_{j} q^{j} \sum_{i}(-1)^{i} dim({\cal C}^{ij})$$ 
$$= \sum_{j} q^{j} \chi({\cal C}^{\bullet,j}) = \sum_{j}q^{j} \chi(H({\cal C}^{\bullet,j})).$$ Here we use the definition of the {\it Euler characteristic of
a chain complex} $$\chi({\cal C}^{\bullet,j}) = \sum_{i}(-1)^{i} dim({\cal C}^{ij})$$ and the fact that the Euler characteristic of the complex is equal to 
the Euler characteristic of its homology. The quantum amplitude associated with the operator $U$ is given 
in terms of the Euler characteristics of the Khovanov homology of the link $K.$ $$\langle K \rangle = \langle \psi |U |\psi \rangle = \sum_{j}q^{j}
\chi(H({\cal C}^{\bullet,j}(K))).$$
\end{enumerate}

Our reformulation of the bracket polynomial in terms of the unitary operator $U$ leads to a new viewpoint on the Khovanov homology
as a representation space for the action of $U.$ The bracket polynomial is then a quantum amplitude that expressing the Euler characteristics of the 
homology associated with this action. The decomposition of the chain complex into the parts ${\cal C}^{i,j}(K)$ corresponds to the eigenspace
decomposition of the operator $U.$ The reader will note that in this case the operator $U$ is already diagonal in the basis of enhanced states for the chain
complex ${\cal C}(K).$  We regard this reformulation as a guide to further questions about the relationship of the Khovanov homology with
quantum information associated with the link $K.$
\bigbreak

As we shall see in the next section, the internal combinatorial structure of the set of enhanced states for the bracket summation leads to the Khovanov
homology theory, whose graded Euler characteristic yields the bracket state sum. Thus we have a quantum statistical interpretation of the Euler
characteristics of the Khovanov homology theory, and a conceptual puzzle about the nature of this relationship with the Hilbert space of that quantum
theory. It is that relationship that is the subject of this paper. The unusual point about the Hilbert space is the each of its basis elements has a
specific combinatorial structure that is related to the topology of the knot $K.$ Thus this Hilbert space is, from the point of view of its basis elements, a
form of ``taking apart" of the topological structure of the knot that we are interested in studying. 
\bigbreak

\noindent {\bf Homological structure of the unitary transformation.}
We now prove a general result about the structure of a chain complex that is also a finite dimensional Hilbert space.
Let ${\cal C}$ be a chain complex over the complex numbers with boundary operator $$\partial: {\cal C}^{i} \longrightarrow {\cal C}^{i+1},$$
with ${\cal C}$ denoting the direct sum of all the ${\cal C}^{i}$, $i = 0,1,2,\cdots n$ (for some $n$). Let 
$$U:{\cal C} \longrightarrow {\cal C}$$ be a unitary operator that satisfies the equation $U \partial + \partial U = 0.$ We do not assume a second
grading $j$ as occurs in the Khovanov homology. However, since $U$ is unitary, it follows \cite{Lang} that there is a basis for ${\cal C}$ in which
$U$ is diagonal. Let ${\cal B} = \{|s\rangle \}$ denote this basis. Let $\lambda_{s}$ denote the eigenvalue of $U$ corresponding to $|s\rangle$ so
that $U |s\rangle = \lambda_{s} |s\rangle.$ Let $\alpha_{s,s'}$ be the matrix element for $\partial$ so that 
$$\partial |s\rangle = \sum_{s'}\alpha_{s,s'} |s' \rangle$$ where $s'$ runs over a set of basis elements so that $i(s') = i(s) + 1.$
\bigbreak

\noindent {\bf Lemma.} With the above conventions, we have that for $|s' \rangle$ a basis element such that $\alpha_{s,s'} \ne 0$ then
$\lambda_{s'} = - \lambda_{s}.$
\smallbreak

\noindent {\bf Proof.} Note that $$U \partial |s \rangle = U(\sum_{s'} \alpha_{s,s'} |s' \rangle) = \sum_{s'} \alpha_{s,s'} \lambda_{s'}|s' \rangle$$ 
while $$\partial U |s \rangle = \partial \lambda_{s} |s \rangle = \sum_{s'} \alpha_{s,s'} \lambda_{s}|s' \rangle.$$
Since $U \partial + \partial U = 0,$ the conclusion of the Lemma follows from the independence of the elements in the basis for the Hilbert space. //
\bigbreak

\noindent In  this way we see that eigenvalues will propagate forward from ${\cal C}^{0}$ with alternating signs according to the appearance
of successive basis elements in the boundary formulas for the chain complex. Various states of affairs are possible in general, with new eigenvaluues
starting at some ${\cal C}^{k}$ for $k > 0.$ The simplest state of affairs would be if all the possible eigenvalues (up to multiplication by $-1$) for $U$
occurred in ${\cal C}^{0}$ so that  $${\cal C}^{0} =\oplus_{\lambda}{\cal C}^{0}_{\lambda}$$ where $\lambda$ runs over all the distinct eigenvalues of $U$
restricted to ${\cal C}^{0},$
and ${\cal C}^{0}_{\lambda}$ is spanned by all $|s \rangle$ in ${\cal C}^{0}$ with $U |s \rangle = \lambda |s \rangle.$ Let us take the further  
assumption that for each $\lambda$ as above, the subcomplexes
$${\cal C}^{\bullet}_{\lambda}: {\cal C}^{0}_{\lambda} \longrightarrow {\cal C}^{1}_{-\lambda} \longrightarrow {\cal C}^{2}_{+\lambda} \longrightarrow \cdots
{\cal C}^{n}_{(-1)^{n}\lambda}$$ have ${\cal C} = \oplus_{\lambda} {\cal C}^{\bullet}_{\lambda} $ as their direct sum. With this assumption about the chain
complex, define 
$|\psi \rangle = \sum_{s} |s\rangle$ as before, with $| s \rangle$ running over the whole basis for ${\cal C}.$  Then it follows just as in the beginning of
this section that $$\langle \psi|U| \psi \rangle = \sum_{\lambda} \lambda \chi(H(C^{\bullet}_{\lambda})).$$ Here $\chi$ denotes the Euler characteristic
of the homology. The point is, that this formula for $\langle \psi|U| \psi \rangle$ takes exactly the form we had for the special case of Khovanov
homology (with eigenvalues $(-1)^{i}q^{j}$), but here the formula occurs just in terms of the eigenspace decomposition of the unitary transformation
$U$ in relation to the chain complex. Clearly there is more work to be done here and we will return to it in a subsequent paper.
\bigbreak 
 
\noindent {\bf Remark on the density matrix.}
Given the state $| \psi \rangle,$ we can define the {\it density matrix} $$\rho = | \psi \rangle \langle \psi|.$$ With this definition it is immediate
that $$Tr(U \rho) = \langle \psi|U| \psi \rangle$$  where $Tr(M)$ denotes the trace of a matrix $M.$ Thus we can restate the form of our result about
Euler characteristics as $$Tr(U \rho) = \sum_{\lambda} \lambda \chi(H(C^{\bullet}_{\lambda})).$$ In searching for an interpretation of the Khovanov
complex in this quantum context it is useful to use this reformulation. For the bracket we have
$$\langle K \rangle = <\psi|U|\psi> = Tr(U \rho).$$
\bigbreak

\noindent{\bf Remark on physics.}
We are taking the whole Hilbert space as the state space of a physical system. If the system were a classical
system then the energetic states of the classical system would correspond to the basis we have chosen for the Hilbert space. This is in line with the 
concepts of quantum mechanics, since the classical states are then in correspondence with a basis for measurement. In quantum context we think of the
elements in the Hilbert space as corresponding to superpositions of possible measurements. All this is then transposed to topological configurations for the
knot or link. But the new ingredient is the one from topology that will make the states of the knot into a chain complex and lead to homological
computations.  Since we are using enhanced states,
each state is a generator for the chain complex. Thus we can regard $\cal{C}(K)$ itself as the chain complex (over the complex numbers) and add in extra
grading structure to boot. There is also the fine structure of the underlying state circles, but this is not seen by the Hilbert space itself and the chain
complex structure can certainly be written just in terms of the basis of the Hilbert space. We want to know what is the
relationship between the unitary transformation and the homological structure. This relation is given by the way the grading works in relation
to $(-1)^{i}q^{j}$ for each  state: the constancy of the quantum grading $j$ under the action of the differential. It is this that makes $ \langle \psi
|U|\psi\rangle$ a graded Euler characteristic for the homology. In this way the unitary transformation is linked with the structure of state transitions that
govern the homology.  At the end of this section we indicate a more general approach to this pattern.
\bigbreak

The key conceptual issue is the construction of a Hilbert space whose basis is the set of observable states of a physical system.
We can always do this for a statistical mechanical system. Let $s$ denote such states and $E(s)$ denote the energy of a state $s$ (observable and real).
Define a unitary transformation $U(t)$ (where $t$ denotes time) by $$U(t)|s \rangle = e^{(it/\hbar)E(s)} |s \rangle.$$ Define the amplitude
$$A(t) = \Sigma_{s} e^{(it/\hbar)E(s)}.$$ Note that we have the formula $$A(t) = \langle \psi|U(t) | \psi \rangle,$$ where 
$|\psi \rangle = \Sigma_{s} |s \rangle$ just as before. Now $U(t)$ can be regarded as the quantum evolution of the state of the system and we have written
this amplitude in analogy with a partition function for the statistical mechanics model.
\bigbreak

One more conceptual problem: We consider an amplitude of the form 
$$\langle \psi | U | \psi \rangle = \sum_{s} (-1)^{i(s)} q^{j(s)}$$ where $q = e^{i \theta}$ is a unit complex number. The question is, how does this fit the
pattern
$$A(t) = \sum_{s} e^{(it/\hbar)E(s)}?$$ So we want $$(t/\hbar)E(s) = \pi i(s) + \theta j(s)$$ which gives 
$$E(s) = (\hbar/t)(\pi i(s) + \theta j(s)).$$ Evolution from time $t=0$ is out of the question. Probably the best interpretation is to think of the
evolution starting at $t=1.$
\bigbreak

In this point of view, the Hilbert space for expressing the bracket polynomial as a quantum statistical amplitude is quite naturally the chain complex for
Khovanov homology with complex coefficients, and the unitary transformation that is the structure of the bracket polynomial acts on the homology of this 
chain complex. This means that the homology classes contain information preserved by the quantum process that underlies the bracket polynomial.
We would like to exploit this direct relationship between the quantum model and the Khovanov homology to obtain deeper information about the relationship of
topology and quantum information theory, and we would like to use this relationship to probe the properties of these topological invariants.
\bigbreak

\section{Background on Khovanov Homology}

In this section, we describe Khovanov homology
along the lines of \cite{Kh,BN}, and we tell the story so that the gradings and the structure of the differential emerge in a natural way.
This approach to motivating the Khovanov homology uses elements of Khovanov's original approach, Viro's use of enhanced states for the bracket
polynomial \cite{Viro}, and Bar-Natan's emphasis on tangle cobordisms \cite{BN2}. We use similar considerations in our paper \cite{DKM}.
\bigbreak

Two key motivating ideas are involved in finding the Khovanov invariant. First of all, one would like to {\it categorify} a link polynomial such as
$\langle K \rangle.$ There are many meanings to the term categorify, but here the quest is to find a way to express the link polynomial
as a {\it graded Euler characteristic} $\langle K \rangle = \chi_{q} \langle {\cal H}(K) \rangle$ for some homology theory associated with $\langle K \rangle.$
\bigbreak

We will use the bracket polynomial and its enhanced states as described in the previous sections of this paper.
To see how the Khovanov grading arises, consider the form of the expansion of this version of the 
bracket polynomial in enhanced states. We have the formula as a sum over enhanced states $s:$
$$\langle K \rangle = \sum_{s} (-1)^{i(s)} q^{j(s)} $$
where $i(s)$ is the number of $B$-type smoothings in $s$, $\lambda(s)$ is the number of loops in $s$ labeled $1$ minus the number of loops
labeled $X,$ and $j(s) = i(s) + \lambda(s)$.
This can be rewritten in the following form:
$$\langle K \rangle  =  \sum_{i \,,j} (-1)^{i} q^{j} dim({\cal C}^{ij}) $$
where we define ${\cal C}^{ij}$ to be the linear span (over the complex numbers for the purpose of this paper, but over the integers or the integers modulo
two for other contexts) of the set of enhanced states with
$i(s) = i$ and $j(s) = j.$ Then the number of such states is the dimension $dim({\cal C}^{ij}).$ 
\bigbreak

\noindent We would like to have a  bigraded complex composed of the ${\cal C}^{ij}$ with a
differential
$$\partial:{\cal C}^{ij} \longrightarrow {\cal C}^{i+1 \, j}.$$ 
The differential should increase the {\it homological grading} $i$ by $1$ and preserve the 
{\it quantum grading} $j.$
Then we could write
$$\langle K \rangle = \sum_{j} q^{j} \sum_{i} (-1)^{i} dim({\cal C}^{ij}) = \sum_{j} q^{j} \chi({\cal C}^{\bullet \, j}),$$
where $\chi({\cal C}^{\bullet \, j})$ is the Euler characteristic of the subcomplex ${\cal C}^{\bullet \, j}$ for a fixed value of $j.$
\bigbreak

\noindent This formula would constitute a categorification of the bracket polynomial. Below, we
shall see how {\it the original Khovanov differential $\partial$ is uniquely determined by the restriction that $j(\partial s) = j(s)$ for each enhanced state
$s$.} Since $j$ is 
preserved by the differential, these subcomplexes ${\cal C}^{\bullet \, j}$ have their own Euler characteristics and homology. We have
$$\chi(H({\cal C}^{\bullet \, j})) = \chi({\cal C}^{\bullet \, j}) $$ where $H({\cal C}^{\bullet \, j})$ denotes the homology of the complex 
${\cal C}^{\bullet \, j}$. We can write
$$\langle K \rangle = \sum_{j} q^{j} \chi(H({\cal C}^{\bullet \, j})).$$
The last formula expresses the bracket polynomial as a {\it graded Euler characteristic} of a homology theory associated with the enhanced states
of the bracket state summation. This is the categorification of the bracket polynomial. Khovanov proves that this homology theory is an invariant
of knots and links (via the Reidemeister moves of Figure 1), creating a new and stronger invariant than the original Jones polynomial.
\bigbreak

\begin{figure}
     \begin{center}
     \begin{tabular}{c}
     \includegraphics[width=6cm]{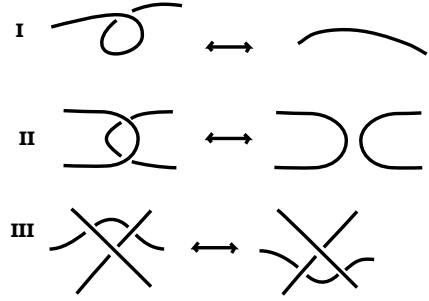}
     \end{tabular}
     \caption{\bf Reidemeister Moves}
     \label{Figure 1}
\end{center}
\end{figure}

We will construct the differential in this complex first for mod-$2$ coefficients. 
The differential is based on regarding two states as {\it adjacent} if one differs from the other by a single smoothing at some site.
Thus if $(s,\tau)$ denotes a pair consisting in an enhanced state $s$ and site $\tau$ of that state with $\tau$ of type $A$, then we consider
all enhanced states $s'$ obtained from $s$ by smoothing at $\tau$ and relabeling only those loops that are affected by the resmoothing.
Call this set of enhanced states $S'[s,\tau].$ Then we shall define the {\it partial differential} $\partial_{\tau}(s)$ as a sum over certain elements in
$S'[s,\tau],$ and the differential by the formula $$\partial(s) = \sum_{\tau} \partial_{\tau}(s)$$ with the sum over all type $A$ sites $\tau$ in $s.$
It then remains to see what are the possibilities for $\partial_{\tau}(s)$ so that $j(s)$ is preserved.
\bigbreak

\noindent Note that if $s' \in S'[s,\tau]$, then $i(s') = i(s) + 1.$ Thus $$j(s') = i(s') + \lambda(s') = 1 + i(s) + \lambda(s').$$ From this
we conclude that $j(s) = j(s')$ if and only if $\lambda(s') = \lambda(s) - 1.$ Recall that 
$$\lambda(s) = [s:+] - [s:-]$$ where $[s:+]$ is the number of loops in $s$ labeled $+1,$ $[s:-]$ is the number of loops
labeled $-1$ (same as labeled with $X$) and $j(s) = i(s) + \lambda(s)$.
\bigbreak

\noindent {\bf Proposition.} The partial differentials $\partial_{\tau}(s)$ are uniquely determined by the condition that $j(s') = j(s)$ for all $s'$
involved in the action of the partial differential on the enhanced state $s.$ This unique form of the partial differential can be described by the 
following structures of multiplication and comultiplication on the algebra \cal{A} = $k[X]/(X^{2})$ where $k = Z/2Z$ for mod-2 coefficients, or $k = Z$
for integral coefficients.
\begin{enumerate}
\item The element $1$ is a multiplicative unit and $X^2 = 0.$
\item $\Delta(1) = 1 \otimes X + X \otimes 1$ and $\Delta(X) = X \otimes X.$
\end{enumerate}
These rules describe the local relabeling process for loops in a state. Multiplication corresponds to the case where two loops merge to a single loop, 
while comultiplication corresponds to the case where one loop bifurcates into two loops.
\bigbreak

\noindent {\bf Proof.}
Using the above description of the differential, suppose that
there are two loops at $\tau$ that merge in the smoothing. If both loops are labeled $1$ in $s$ then the local contribution to $\lambda(s)$ is $2.$
Let $s'$ denote a smoothing in $S[s,\tau].$ In order for the local $\lambda$ contribution to become $1$, we see that the merged loop must be labeled $1$.
Similarly if the two loops are labeled $1$ and $X,$ then the merged loop must be labeled $X$ so that the local contribution for $\lambda$ goes from 
$0$ to $-1.$ Finally, if the two loops are labeled $X$ and $X,$ then there is no label available for a single loop that will give $-3,$ so we define
$\partial$ to be zero in this case. We can summarize the result by saying that there is a multiplicative structure $m$ such that 
$m(1,1) = 1, m(1,X) = m(X,1) = x, m(X,X) = 0,$ and this multiplication describes the structure of the partial differential when two loops merge.
Since this is the multiplicative structure of the algebra ${\cal A} = k[X]/(X^{2}),$ we take this algebra as summarizing the differential.
\bigbreak

Now consider the case where $s$ has a single loop at the site $\tau.$ Smoothing produces two loops. If the single loop is labeled $X,$ then we must label
each of the two loops by $X$ in order to make $\lambda$ decrease by $1$. If the single loop is labeled $1,$ then we can label the two loops by
$X$ and $1$ in either order. In this second case we take the partial differential of $s$ to be the sum of these two labeled states. This structure
can be described by taking a coproduct structure with $\Delta(X) = X \otimes X$ and $\Delta(1) = 1 \otimes X + X \otimes 1.$
We now have the algebra ${\cal A} = k[X]/(X^{2})$ with product $m: {\cal A} \otimes {\cal A} \longrightarrow {\cal A}$ and coproduct
$\Delta: {\cal A} \longrightarrow {\cal A} \otimes {\cal A},$ describing the differential completely. This completes the proof. //
\bigbreak

Partial differentials are defined on each enhanced state $s$ and a site $\tau$ of type
$A$ in that  state. We consider states obtained from the given state by  smoothing the given site $\tau$. The result of smoothing $\tau$ is to
produce a new state $s'$ with one more site of type $B$ than $s.$ Forming $s'$ from $s$ we either amalgamate two loops to a single loop at $\tau$, or
we divide a loop at $\tau$ into two distinct loops. In the case of amalgamation, the new state $s$ acquires the label on the amalgamated circle that
is the product of the labels on the two circles that are its ancestors in $s$. This case of the partial differential is described by the
multiplication in the algebra. If one circle becomes two circles, then we apply the coproduct. Thus if the circle is labeled $X$, then the resultant
two circles are each labeled $X$ corresponding to $\Delta(X) = X \otimes X$. If the orginal circle is labeled $1$ then we take the partial boundary
to be a sum of two enhanced states with  labels $1$ and $X$ in one case, and labels $X$ and $1$ in the other case,  on the respective circles. This
corresponds to $\Delta(1) = 1 \otimes X + X \otimes 1.$ Modulo two, the boundary of an enhanced state is the sum, over all sites of type $A$ in the
state, of the partial boundaries at these sites. It is not hard to verify directly that the square of the  boundary mapping is zero (this is the identity of 
mixed partials!) and that it behaves
as advertised, keeping $j(s)$ constant. There is more to say about the nature of this construction with respect to Frobenius algebras and tangle
cobordisms. In Figures 2,3 and 4 we illustrate how the partial boundaries can be conceptualized in terms of surface cobordisms. The equality of mixed
partials corresponds to topological equivalence of the corresponding surface cobordisms, and to the relationships between Frobenius algebras and the
surface cobordism category. In particular, in Figure 4 we show how in a key case of two sites (labeled 1 and 2 in that Figure) the two orders of partial 
boundary are $$\partial_{2} \partial_{1} = (1 \otimes m) \circ (\Delta \otimes 1)$$ and 
$$\partial_{1} \partial_{2} = \Delta \circ m.$$ In the Frobenius algebra ${\cal A} = k[X]/(X^{2})$ we have the identity
$$(1 \otimes m) \circ (\Delta \otimes 1) = \Delta \circ m.$$ Thus the Frobenius algebra implies the identity of the mixed partials.
Furthermore, in Figure 3 we see that this identity corresponds to the topological equivalence of cobordisms under an exchange of saddle points.
There is more to say about all of this, but we will stop here.
The proof of invariance of Khovanov homology with respect to the Reidemeister moves (respecting grading changes) will not be
given here. See \cite{Kh,BN,BN2}. It is remarkable that this version of Khovanov homology is uniquely specified by natural ideas about adjacency of states
in the bracket polynomial.
\bigbreak

\begin{figure}
     \begin{center}
     \begin{tabular}{c}
     \includegraphics[width=6cm]{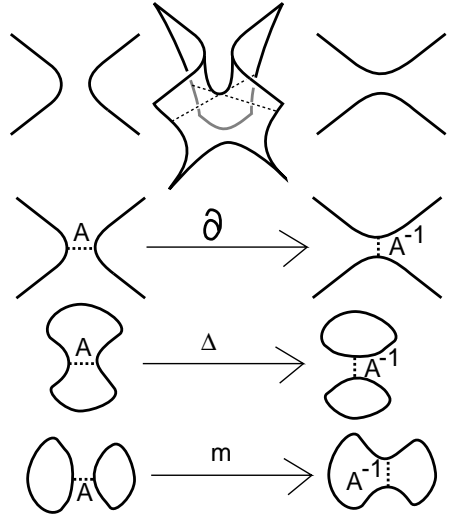}
     \end{tabular}
     \caption{\bf SaddlePoints and State Smoothings}
     \label{Figure 2}
\end{center}
\end{figure}

\begin{figure}
     \begin{center}
     \begin{tabular}{c}
     \includegraphics[width=7cm]{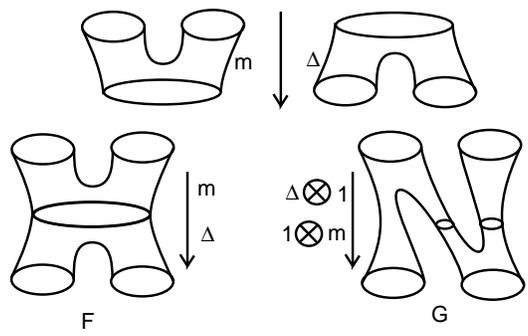}
     \end{tabular}
     \caption{\bf Surface Cobordisms}
     \label{Figure 3}
\end{center}
\end{figure}

\begin{figure}
     \begin{center}
     \begin{tabular}{c}
     \includegraphics[width=9cm]{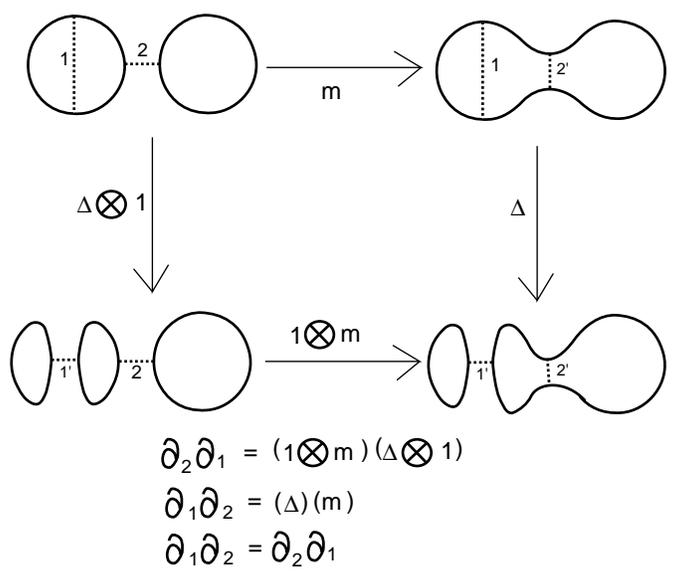}
     \end{tabular}
     \caption{\bf Local Boundaries Commute}
     \label{Figure 4}
\end{center}
\end{figure}

\noindent {\bf Remark on Integral Differentials}. Choose an ordering for the crossings in the link diagram $K$ and denote them by $1,2,\cdots n.$
Let $s$ be any enhanced state of $K$ and let $\partial_{i}(s)$ denote the chain obtained from $s$ by applying a partial boundary at the $i$-th site
of $s.$ If the $i$-th site is a smoothing of type $A^{-1}$, then $\partial_{i}(s) = 0.$ If the $i$-th site is a smoothing of type $A$, then
$\partial_{i}(s)$ is given by the rules discussed above (with the same signs). The compatibility conditions that we have discussed show that
partials commute in the sense that $\partial_{i} (\partial_{j} (s)) = \partial_{j} (\partial_{i} (s))$ for all $i$ and $j.$ One then defines
signed boundary formulas in the usual way of algebraic topology. One way to think of this regards the complex as the analogue of 
a complex in de Rahm cohomology. Let $\{dx_{1}, dx_{2},\cdots, dx_{n}\}$ be a formal basis for a Grassmann algebra so that 
$dx_{i} \wedge dx_{j} = - dx_{j} \wedge dx_{i}$ Starting with enhanced states $s$ in $C^{0}(K)$ (that is, state with all $A$-type smoothings)
Define formally, $d_{i}(s) = \partial_{i}(s) dx_{i}$ and regard $d_{i}(s)$ as identical with $\partial_{i} (s)$ as we have previously regarded it
in $C^{1}(K).$ In general, given an enhanced state $s$ in $C^{k}(K)$ with $B$-smoothings at locations $i_{1} < i_{2} < \cdots < i_{k},$ we represent
this chain as $s \, dx_{i_{1}} \wedge \cdots \wedge dx_{i_{k}}$ and define
$$\partial ( s \, dx_{i_{1}} \wedge \cdots \wedge dx_{i_{k}} ) = \sum_{j=1}^{n} \partial_{j}(s) \, dx_{j} \wedge dx_{i_{1}} \wedge \cdots \wedge dx_{i_{k}},$$
just as in a de Rahm complex. The Grassmann algebra automatically computes the correct signs in the chain complex, and this boundary formula gives the 
original boundary formula when we take coefficients modulo two. Note, that in this formalism, partial differentials $\partial_{i}$ of enhanced states with a
$B$-smoothing at the site $i$ are zero due to the fact that $dx_{i} \wedge dx_{i} = 0$ in the Grassmann algebra. There is more to discuss
about the use of Grassmann algebra in this context. For example, this approach clarifies parts of the construction in \cite{M}.
\bigbreak

It of interest to examine this analogy between the Khovanov (co)homology and de Rahm cohomology. In that analogy the enhanced
states correspond to the differentiable functions on a manifold. The Khovanov complex $C^{k}(K)$ is generated by elements of the form
$s \, dx_{i_{1}} \wedge \cdots \wedge dx_{i_{k}}$ where the enhanced state $s$ has $B$-smoothings at exactly the sites $i_{1},\cdots, i_{k}.$
If we were to follow the analogy with de Rahm cohomology literally, we would define a new complex $DR(K)$ where $DR^{k}(K)$ is generated by elements
$s \, dx_{i_{1}} \wedge \cdots \wedge dx_{i_{k}}$ where $s$ is {\it any} enhanced state of the link $K.$ The partial boundaries are defined in the same
way as before and the global boundary formula is just as we have written it above. This gives a {\it new} chain complex associated with the link $K.$
Whether its homology contains new topological information about the link $K$ will be the subject of a subsequent paper.
\bigbreak

\noindent {\bf A further remark on de Rham cohomology.} There is another deep relation with the de
Rham complex: In \cite{Pr} it was observed that Khovanov homology is related to Hochschild
homology and Hochschild homology is thought to be an algebraic version of de Rham chain
complex (cyclic cohomology corresponds to de Rham cohomology), compare \cite{Lo}.
\bigbreak 

\section {Other Homological States}

The formalism that we have pursued in this paper to relate Khovanov homology and the Jones polynomial with quantum statistics can be generalized
to apply to other situations. It is possible to associate a finite dimensional Hilbert space that is also a chain complex (or cochain complex) to
topological structures other than knots and links. To formulate this, let $X$  denote the topological structure and ${\cal C}(X)$ denote the associated
linear space endowed with a boundary operator $\partial: {\cal C}(X) \longrightarrow {\cal C}(X),$ with $\partial \partial = 0.$ We want to consider 
situations where there is a unitary operator $U:{\cal C}(X) \longrightarrow {\cal C}(X)$ such that $U \partial + \partial U = 0$ so that $U$
induces a unitary action on ${\cal H}(X),$ the homology of ${\cal C}(X)$ with respect to $\partial.$
\bigbreak

For example, let $M$ be a differentiable manifold and ${\cal C}(M)$ denote the DeRham complex of $M$ over the complex numbers. Then for a differential
form of the type $f(x) \omega$ in local coordinates $x_{1}, \cdots ,x_{n}$ and $\omega$ a wedge product of a subset of $dx_{1} \cdots dx_{n},$ we have 
$$d(f \omega) = \sum_{i=1}^{n} (\partial f/\partial x_{i}) dx_{i} \wedge \omega.$$ Here $d$ is the differential for the DeRham complex. Then 
${\cal C}(M)$ has basis the set of $|f(x) \omega \rangle$ where $\omega = dx_{i_{1}} \wedge \cdots \wedge dx_{i_{k}}$ with 
$i_{1} < \cdots < i_{k} .$ We could achieve $U d + d U = 0$ if $U$ is a very simple unitary operator (e.g. multiplication by phases that do not depend on 
the coordinates $x_{i}$) but in general it will be an interesting problem to determine all unitary operators $U$ with this property.
\bigbreak

Even in the case of Khovanov homology, we could keep the homology theory fixed and ask for other unitary operators $U$ that satisfy
$U \partial + \partial U = 0.$ Knowing other examples of such operators would shed light on the nature of Khovanov homology from the point of 
view of quantum statistics.
\bigbreak

\section {Appendix - The Hadamard Test}

In order to make a quantum computation of the trace of a unitary matrix $U$, one can use the {\it Hadamard test} to obtain the diagonal matrix
elements $\langle \psi|U|\psi \rangle$ of $U.$ The trace is then the sum of these matrix elements as $|\psi \rangle$ runs over an orthonormal basis for 
the vector space. In the application to the algorithm described here for the Jones polynomial it is only necessary to compute one number of the form 
$\langle \psi|U|\psi \rangle.$ The Hadamard test proceeds as follows.
\bigbreak

\noindent We first obtain
$$\frac{1}{2} + \frac{1}{2}Re\langle \psi|U|\psi \rangle$$ 
as an expectation by applying the Hadamard gate $H$
$$H|0 \rangle = \frac{1}{\sqrt{2}}(|0\rangle + |1\rangle)$$
$$H|1 \rangle = \frac{1}{\sqrt{2}}(|0\rangle - |1\rangle)$$
to the first qubit of 
$$C_{U} \circ (H \otimes 1) |0 \rangle |\psi \rangle = \frac{1}{\sqrt{2}}(|0\rangle \otimes|\psi \rangle + |1\rangle \otimes U|\psi\rangle.$$
Here $C_{U}$ denotes controlled $U,$ acting as $U$ when the control bit is $|1 \rangle$ and the identity mapping when the control bit is $|0 \rangle.$ We
measure the expectation for the first qubit $|0 \rangle$ of the resulting state
$$(H \otimes 1) \circ C_{U} \circ (H \otimes 1) |0 \rangle |\psi \rangle =$$
$$\frac{1}{2}(H|0\rangle \otimes|\psi \rangle + H|1\rangle \otimes U|\psi\rangle)
=\frac{1}{2}((|0\rangle + |1\rangle) \otimes|\psi \rangle + (|0\rangle - |1\rangle) \otimes U|\psi\rangle)$$
$$=\frac{1}{2}(|0\rangle \otimes (|\psi \rangle + U|\psi\rangle) + |1\rangle \otimes(|\psi \rangle - U|\psi\rangle)).$$
This expectation is $$\frac{1}{4}(\langle \psi | + \langle \psi| U^{\dagger})(|\psi \rangle + U|\psi\rangle) = \frac{1}{2} + \frac{1}{2}Re\langle \psi|U|\psi
\rangle.$$ In Figure 5, we illustrate this computation with a diagram that indicates the structure of the test with parallel lines corresponding to 
tensor products of the single qubit space (with three lines chosen for illustration as the size of $U$). The extra tensor factor is indicated on the top
line with the Hadamard matrix $H$ indicated by a box and the control of $U$ indicated by a circle with a vertical control line extending down to the 
$U$-box. The half -circle on the top line on the right stands for the measurement of that line that is used for the computation. Thus Figure 5 represents
a circuit diagram for the quantum computation of $\frac{1}{2} + \frac{1}{2}Re\langle \psi|U|\psi \rangle$ and hence the quantum computation of 
 $\frac{1}{2} + \frac{1}{2}Re\langle K \rangle$ when $U$ is taken to be the unitary tranformation corresponding to the bracket polynomial, as
discussed in previous sections  of this paper. 
\bigbreak

\noindent The imaginary
part is  obtained by applying the same procedure to 
$$\frac{1}{\sqrt{2}}(|0\rangle \otimes|\psi \rangle - i|1\rangle \otimes U|\psi\rangle$$
This is the method used in
\cite{AJL,QCJP,3Strand}, and the reader may wish to contemplate its efficiency in the context of this simple model. Note that the Hadamard test enables this
quantum  computation to estimate the trace of any unitary matrix $U$ by repeated trials that estimate individual matrix entries $\langle \psi|U|\psi\rangle.$
\bigbreak

\begin{figure}
     \begin{center}
     \begin{tabular}{c}
     \includegraphics[width=9cm]{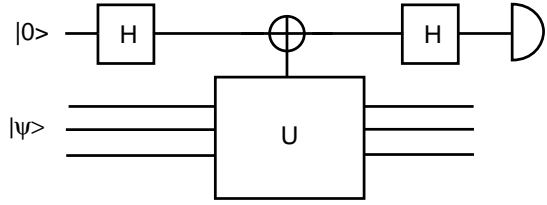}
     \end{tabular}
     \caption{\bf Hadamard Test}
     \label{Figure 5}
\end{center}
\end{figure}

\end{document}